\documentclass[12pt]{article}

\usepackage{amsthm,amsmath,amssymb,url,cite}
\newtheorem{thm}{Theorem}[section]
\newtheorem{lem}[thm]{Lemma}

\theoremstyle{remark}

\newtheorem*{eg*}{Example}

\numberwithin{equation}{section}

\renewcommand{\qed}{{\hfill\rule{4pt}{7pt}}\medskip}

\def\pf{\noindent {\it Proof.} }

\begin{document}

\begin{center}
{\Large\bf Contiguous relations and summation and transformation
formulas for basic hypergeometric series}
\end{center}

\vskip 2mm \centerline{\large Feng Gao$^1$ and Victor J. W.
Guo$^2$\footnote{Corresponding author.} }

\vskip 3mm \centerline{\footnotesize Department of Mathematics, East
China Normal University, Shanghai 200062,
 People's Republic of China}

\vskip 1mm \centerline {\footnotesize {$^1$\tt
 {\tt nature1314@163.com},\qquad  $^2$jwguo@math.ecnu.edu.cn,\quad
http://math.ecnu.edu.cn/\textasciitilde{jwguo}}
 }

\vskip 0.7cm {\small \noindent{\bf Abstract.} By using contiguous relations for basic hypergeometric series,
we give simple proofs of Bailey's $_4\phi_3$ summation, Carlitz's $_5\phi_4$ summation,
Sears' $_3\phi_2$ to $_5\phi_4$ transformation, Sears' ${}_4\phi_3$ transformations, Chen's bibasic summation,
Gasper's split poised $_{10}\phi_9$ transformation, Chu's bibasic symmetric transformation.
Along the same line, finite forms of Sylvester's identity, Jacobi's triple product identity,
and Kang's identity are also obtained.


\vskip 2mm
\noindent{\it MR Subject Classifications}: 33D15}

\section{Introduction}
In the previous work \cite{Guo,Guo2,GZ07}, many terminating summation and transformation formulas for
basic hypergeometric series are proved by using contiguous relations and mathematical induction.
The present paper is a complement to \cite{Guo,Guo2,GZ07}. The difference here is that we will
usually apply mathematical induction twice rather than once.
Recall that the {\it $q$-shifted factorial} is defined by
\[
(a;q)_\infty=\prod_{k=0}^{\infty}(1-aq^{k}),\quad\text{and}
\quad (a;q)_n=\frac{(a;q)_\infty}{(aq^n;q)_\infty} \quad\text{for}
\quad n\in\mathbb{Z}.
\]
We employ the abbreviated notation
$$
(a_1,a_2,\ldots,a_m;q)_n
=(a_1;q)_n (a_2;q)_n \cdots(a_m;q)_n,
\quad\text{for}\quad n=\infty\quad\text{or}\quad n\in\mathbb{Z}.
$$
The {\it basic hypergeometric series $_{r+1}\phi_r$} is defined as
$$
_{r+1}\phi_{r}\left[\begin{array}{c}
a_1,a_2,\ldots,a_{r+1}\\
b_1,b_2,\ldots,b_{r}
\end{array};q,\, z
\right]
=\sum_{k=0}^{\infty}\frac{(a_1,a_2,\ldots,a_{r+1};q)_k z^k}
{(q,b_1,b_2,\ldots,b_{r};q)_k}.
$$
An $_{r+1}\phi_{r}$ series is called {\it well-poised} if
$a_1q=a_2b_1=\cdots=a_{r+1}b_r$.

Let
\begin{align*}
F_k(a_1,a_2,\ldots,a_{r+1};q,z):=\frac{(a_1,\ldots,a_{r+1};q)_k}
{(q,a_1q/a_2,\ldots,a_1q/a_{r+1};q)_k} z^k
\end{align*}
be the $k$-th term of a well-poised hypergeometric series. In \cite{GZ07},
the following contiguous relations were established:
\begin{align}
F_k(a_1,a_2,\ldots,a_rq, a_{r+1};q,z)&-F_k(a_1,a_2,\ldots,a_r, a_{r+1}q;q,z)\nonumber \\
&=\alpha F_{k-1}(a_1q^2,a_2q,\ldots,a_{r+1}q;q,z),
\label{eq:contiguous1} \\[5pt]
F_k(a_1,a_2,\ldots,a_r, a_{r+1};q,qz)&-F_k(a_1,a_2,\ldots,a_r, a_{r+1}q;q,z)\nonumber\\
&=\beta F_{k-1}(a_1q^2,a_2q,\ldots,a_{r+1}q;q,z), \label{eq:contiguous2}
\end{align}
where
\begin{align*}
\alpha &=\frac{(a_r-a_{r+1})(1-a_1/a_ra_{r+1})(1-a_1)(1-a_1q)(1-a_2)
\cdots (1-a_{r-1})z}{(1-a_1/a_r)(1-a_1/a_{r+1})(1-a_1q/a_2)
\cdots (1-a_1q/a_{r+1})}, \\[5pt]
\beta &=-\frac{(1-a_1)(1-a_1q)(1-a_2)\cdots (1-a_r)z}
{(1-a_1/a_{r+1})(1-a_1q/a_2)\cdots (1-a_1q/a_{r+1})},
\end{align*}
and quite a few terminating summation and transformation formulas for basic hypergeometric series
were proved by applying \eqref{eq:contiguous1} or \eqref{eq:contiguous2}.
Note that Krattenthaler \cite{Krattenthaler0} has indicated how to derive contiguous relations
from special cases of known formulas (available at http://www.mat.univie.ac.at/\textasciitilde kratt/papers.html).

In this paper we shall prove more such identities based on the contiguous relations
\eqref{eq:contiguous1}, \eqref{eq:contiguous2} or others. Suppose that
\begin{align}
\sum_{k=0}^{n}F_{n,k}(a_1,\ldots,a_s)=S_n(a_1,\ldots,a_s),
\label{eq:fnk-sum}
\end{align}
where $F_{n,k}(a_1,\ldots,a_s)=0$ if $k<0$ or $k>n$. If one can show
that the summand $F_{n,k}(a_1,\ldots,a_s)$ satisfies the following
recurrence relation:
\begin{align}
F_{n,k}(a_1q,a_2\ldots,a_s)-F_{n,k}(a_1,a_2,\ldots,a_s) =\gamma_n
F_{n-v,k-1}(b_1,\ldots,b_s),\ \text{$v=1$ or $2$} \label{eq:fnk-ind}
\end{align}
for some parameters $b_1,\ldots,b_s$, where $\gamma_n $ is
independent of $k$, then we can prove \eqref{eq:fnk-sum} by
induction on $n$. Note that \eqref{eq:ssn} is not a special
case of Sister Celine's method \cite[p. 58, (4.3.1)]{PWZ}, as mentioned by \cite{GZ07}.

Of course, we need to check that $S_n(a_1,\ldots,a_s)$ satisfies the
following recurrence relation
\begin{align}
S_{n}(a_1q,a_2,\ldots,a_s)-S_{n}(a_1,a_2,\ldots,a_s)=\gamma_n
S_{n-v}(b_1,\ldots,b_s). \label{eq:ssn}
\end{align}

If $S_n(a_1,\ldots,a_s)$ appears as a {\it closed form} as in
Bailey's $_4\phi_3$ summation formula \eqref{eq:cr002},
 then the verification of \eqref{eq:ssn} is quite easy.
If $S_n(a_1,\ldots,a_s)$ is of the following form:
\begin{align*}
S_n(a_1,\ldots,a_s)=\sum_{k=0}^{n}G_{n,k}(a_1,\ldots,a_s),
%\label{eq:gnk-sum}
\end{align*}
where $G_{n,k}(a_1,\ldots,a_s)=0$ if $k<0$ or $k>n$, then we may try to apply
$q$-Gosper's algorithm \cite[p.~75]{Koepf} to find a sequence
$H_{n,k}(a_1,\ldots,a_s)$ of closed forms such that
\begin{align}
%&\hskip -3mm
G_{n,k}(a_1q,a_2,\ldots,a_s)-G_{n,k}(a_1,a_2\ldots,a_s)-\gamma_n
G_{n-v,k-1}(b_1,\ldots,b_s)
%\nonumber\\[5pt]
=H_{n,k}-H_{n,k-1},\label{eq:hnkk}
\end{align}
where $H_{n,n}=H_{n,-1}=0$. If \eqref{eq:hnkk} exists, then by telescoping we get \eqref{eq:ssn}.

It is easy to see that almost all identities are trivial for $n=0$ or $n=1$.
For the induction step, we need firstly to prove that \eqref{eq:fnk-sum} is true for some special $a_1$, and
secondly give \eqref{eq:fnk-ind} (or, in addition, \eqref{eq:hnkk}). In this way, we shall give simple proofs
of Bailey's $_4\phi_3$ summation, Carlitz's $_5\phi_4$ summation,
Sears' $_3\phi_2$ to $_5\phi_4$ transformation, Sears' ${}_4\phi_3$ transformations, Chen's bibasic summation,
Gasper's split poised $_{10}\phi_9$ transformation, Chu's bibasic symmetric transformation.
Moreover, we will give a finite form of Sylvester's identity, as well as a finite form of
Kang's identity \cite{Kang}.

Note that all terminating identities can be also proved automatically by
using the $q$-Zeilberger algorithm (see, for example, \cite{BK,Koepf,PWZ}).
Moreover, the so-called ``Abel's method" (see, for example, \cite{Chu,CJ})
is also a nice method to deal with such identities. Gasper \cite{Gasper}
has already developed elementary proofs for many summation formulas for basic hypergeometric series
by considering contiguous relations.  However, our method is a little different from Gasper's method,
since the latter is mainly applied to prove nonterminating basic hypergeometric series, while
our method is usually applied to prove terminating series.

%%%%%%%%%%%%%%%%%%%%%%%%%%%%%%%%%%%%%%%%%%%%%%%%%%%%%%%%%%%%%%%%%%%%%%%%%%%%%%%%%
\section{Bailey's $_4\phi_3$ summation formula}

Bailey \cite{Bailey41} obtained the following $_4\phi_3$ summation formula, which was later generalized
by Carlitz \cite{Carlitz}.

\begin{thm}[Bailey's $_4\phi_3$ summation formula]\label{thm:cr002}For $n\geq 0$, there holds
\begin{align}
 _4\phi_3
\left[\begin{array}{c}
q^{-n},\, b,\, c,\, -q^{1-n}/bc \\
q^{1-n}/b,\, q^{1-n}/c,\, -bc\end{array};q,\,q \right]
 &=\begin{cases}
\displaystyle \frac{(q, b^2, c^2;q^2)_m(bc;q)_{2m}}{(b,c;q)_{2m}(b^2c^2;q^2)_m},
&\text{if $n = 2m$,}\\[10pt]
0, &\text{if $n = 2m+1$.} \label{eq:cr002}
\end{cases}
\end{align}
\end{thm}

The following lemma is the $c=q$ case of Theorem \ref{thm:cr002}.
\begin{lem}\label{lem:c=q}
For $n\geq 0$, there holds
\begin{align}
&\sum_{k=0}^{n}(-1)^k\frac{(b;q)_k(b;q)_{n-k}}{(-b;q)_{n-k+1} (-bq;q)_k}
=
\begin{cases}
\displaystyle\frac{(b;q)_{n+1}}{(-bq;q)_{n}(1-b^2q^{n})}, &\text{if $n$ is even,} \\[10pt]
0, &\text{if $n$ is odd.}
\end{cases} \label{eq:case2}
\end{align}
\end{lem}

\pf Since the $k$-th term and $(n-k)$-th term only differ by a factor $(-1)^{n}$, the left-hand side of
\eqref{eq:case2} is equal to $0$ if $n$ is odd.

Now consider the case where $n$ is even. Since
\begin{align*}
(1-b^2q^{n})=(1-bq^{n-k})-(1-bq^k)+(1+bq^{n-k})(1-bq^k),
\end{align*}
we have
\begin{align}
&\hskip -3mm
(1-b^2q^{n})\sum_{k=0}^{n}(-1)^k\frac{(b;q)_k(b;q)_{n-k}}{(-b;q)_{n-k+1}(-bq;q)_k} \nonumber \\[5pt]
&=\sum_{k=0}^{n}(-1)^k\frac{(b;q)_k(b;q)_{n+1-k}-(b;q)_{k+1}(b;q)_{n-k}}{(-b;q)_{n-k+1}(-bq;q)_k}
+\sum_{k=0}^{n}(-1)^k\frac{(b;q)_{k+1}(b;q)_{n-k}}{(-b;q)_{n-k}(-bq;q)_k}.  \label{eq:new1}
\end{align}
Noticing that in the first summation of \eqref{eq:new1} the $k$-th term and the $(n-k)$-th term only differ by a sign,
we immediately get
\begin{align}
\sum_{k=0}^{n}(-1)^k\frac{(b;q)_k(b;q)_{n+1-k}-(b;q)_{k+1}(b;q)_{n-k}}
{(-b;q)_{n-k+1}(-bq;q)_k}=0.  \label{eq:new2}
\end{align}
For the same reason, we have
\begin{align}
\sum_{k=0}^{n-1}(-1)^k\frac{(b;q)_{k+1}(b;q)_{n-k}}{(-b;q)_{n-k}
(-b;q)_{k+1}}=0. \label{eq:new3}
\end{align}
The proof then follows from combining \eqref{eq:new1}--\eqref{eq:new3}. \qed

\noindent{\it Proof of Theorem \ref{thm:cr002}. } It  suffices to
prove it for the cases $c=q^{M}$ ($M\geq 1$). We proceed by
induction on $M$. For $M=1$, Equation \eqref{eq:cr002} reduces to
\eqref{eq:case2}.

Assume that \eqref{eq:cr002} holds for $c=q^M$. Let
$$
F_{n,k}(b,c,q)
=\frac{(q^{-n};q)_k(b;q)_k(c;q)_k(-q^{1-n}/bc;q)_kq^k}{(q;q)_k(q^{1-n}/b;q)_k(q^{1-n}/c;q)_k(-bc;q)_k}.
$$
Applying the contiguous relation  (1.1) with $
a_r=b $ and $ a_{r+1}=-q^{-n}/bc $, we have
\begin{align}
&F_{n,k}(bq,c,q)-F_{n,k}(b,c,q)=\alpha_n F_{n-2,k-1}(bq,cq,q), \label{eq:fnk-bc}
\end{align}
where
$$
\alpha_n=\frac{(b+q^{-n}/bc)(1-c^2)(1-q^{-n})(1-q^{1-n})q}{(1-q^{-n}/b)(1-q^{1-n}/b)(1-q^{1-n}/c)(1+bc)(1+bcq)}.
$$
Let
$$
S_{n}(b,c,q):=\sum_{k=0}^n F_{n,k}(b,c,q).
$$
Then summing \eqref{eq:fnk-bc} over $k$ from $0$ to $n$ gives
$$
S_{n-2}(bq,cq,q)=\alpha_n^{-1}(S_{n}(bq,c,q)-S_{n}(b,c,q)),
$$
from which one can readily check that
\eqref{eq:cr002} holds for $cq=q^{M+1}$.
\qed

%%%%%%%%%%%%%%%%%%%%%%%%%%%%%%%%%%%%%%%%%%%%%%%%%%%%%%%%%%%%%%%%%%%%%%%%%%%%%%%%%%%%%%%%%%%%%%%%
\section{Carlitz's $_5\phi_4$ summation formula}
By using the $q$-Pfaff-Saalsch\"utz formula,
Carlitz \cite{Carlitz} obtained the following $_5\phi_4$ summation formula, which
is a generalization of \cite[Eq. (1)]{Jackson21}. Simple proofs of Carlitz's formula
have already been given by Guo \cite{Guo} and Chu and Jia \cite{CJ}. Here we give a
simpler one.
\begin{thm}[Carlitz's $_5\phi_4$ summation formula]For $n\geq 0$, there holds
\begin{align}
&\hskip -2mm
 _5\phi_4
\left[\begin{array}{c}
q^{-n},\, b,\, c,\, d,\, e \\
q^{1-n}/b,\, q^{1-n}/c,\, q^{1-n}/d,\, q^{1-n}/e\end{array};q,\,q \right] \nonumber \\
&=q^{m(1+m-n)}(de)^{-m}
\frac{(q^{-m};q)_{2m}(q^{1-n}/bc,q^{1-n}/bd,q^{1-n}/be;q)_{m}(q^{2m-n};q)_{n-2m}}
{(q,q^{1-n}/b,q^{1-n}/c,q^{1-n}/d,q^{1-n}/e,q^{n-m}c;q)_{m}},  \label{eq:32-54}
\end{align}
where \,$bcde=q^{1+m-2n}$ and $m=\lfloor n/2\rfloor$.
\end{thm}

\pf
It suffices to prove the cases $d=-q^{-M}$ ($M\geq \lfloor n/2\rfloor$). We
proceed by induction on $M$. For $M=\lfloor n/2\rfloor$, the identity
\eqref{eq:32-54} reduces to \eqref{eq:cr002}. Assume that \eqref{eq:32-54} holds for $d=-q^{-M}$.

Let
$$
F_{n,k}(b,c,d,q) =\frac{(q^{-n},b,c,d,q^{1-\lfloor
(3n+1)/2\rfloor}/bcd;q)_{k}}
{(q,q^{1-n}/b,q^{1-n}/c,q^{1-n}/d,q^{\lfloor
(n+1)/2\rfloor}bcd;q)_k}q^k.
$$
Applying the contiguous relation
\eqref{eq:contiguous1}, we have
\begin{align}
F_{n,k}(b,c,dq,q)- F_{n,k}(b,c,d,q)=\alpha_n
F_{n-2,k-1}(bq,cq,dq,q),  \label{eq:fnkbcd}
\end{align}
where
\begin{align*}
\alpha_n
=-\frac{d(1-b)(1-c)(1-q^{\lfloor (3n+1)/2\rfloor}bcd^2)(1-q^{\lfloor (n+1)/2\rfloor}bc)(1-q^n)(1-q^{n-1})
q^{\lfloor n/2\rfloor-1}}
{(1-q^{\lfloor (n+1)/2\rfloor}bcd)(1-q^{\lfloor (n+3)/2\rfloor}bcd)(1-q^{n-1}b)(1-q^{n-1}c)(1-q^{n-1}d)(1-q^{n}d)}.
\end{align*}
Let
$$
S_{n}(b,c,d,q):=\sum_{k=0}^n F_{n,k}(b,c,d,q).
$$
Then summing \eqref{eq:fnkbcd} over $k$ from $0$ to $n$ gives
$$
S_{n}(b,c,d,q)=S_{n}(b,c,dq,q)-\alpha_n S_{n-2}(bq,cq,dq,q),
$$
from which one can verify that
\eqref{eq:32-54} holds for $d=-q^{-M-1}$. \qed

%%%%%%%%%%%%%%%%%%%%%%%%%%%%%%%%%%%%%%%%%%%%%%%%%%%%%%%%%%%%%%%%%%%%%%%%%%%%%%%%%%%%%%%%%%%%%%%%%
\section{Sears' $_3\phi_2$ to $_5\phi_4$ transformation}
The following transformation was first obtained by Sears \cite[(4.1)]{Sears}.
See also Carlitz \cite[(2.4)]{Carlitz}.
\begin{thm}[Sears' $_3\phi_2$ to $_5\phi_4$ transformation]\label{thm:1.2}
For $ a=q^{-n}$, $n=0,1,2,\ldots,$ there holds
\begin{align}
_3\phi_2 \left[\begin{array}{c}
a,\, b,\, c \\
aq/b,\, aq/c\end{array}; q,\,\frac{aqx}{bc}
\right]
=\frac{(ax;q)_{\infty}}{(x;q)_{\infty}}  {}_5\phi_4
\left[\begin{array}{c}
a^{\frac{1}{2}}, \, -a^{\frac{1}{2}},\, (aq)^{\frac{1}{2}},\, -(aq)^{\frac{1}{2}},\, aq/bc \\
aq/b,\,aq/c,\, ax,\, q/x,\end{array};q,\,q\right].
\label{eq:32to54}
\end{align}
\end{thm}

\pf
Let
$$
F_{n,k}(x,b,c,q)
=\frac{(q^{-n},b,c;q)_k}{(q^{1-n}/b,q^{1-n}/c,q;q)_k}\left(\frac{q^{1-n}x}{bc}\right)^k.
$$
Applying the contiguous relation \eqref{eq:contiguous2} with $a_{r+1}=c$ and $z=xq^{-n}/bc$,
we have
\begin{align*}
F_{n,k}(x,b,c,q)-F_{n,k}(x,b,cq,q)=\beta_n F_{n-2,k-1}(xq^{-1},bq,cq,q),
\end{align*}
where
$$
\beta_n=-\frac{(1-q^{-n})(1-q^{1-n})(1-b)xq^{-n}/bc}{(1-q^{-n}/c)(1-q^{1-n}/b)(1-q^{1-n}/c)}.
$$

Let
\begin{align*}
G_{n,k}(x,b,c,q)
&=\frac{(q^{-n}x;q)_n(q^{-n};q)_{2k}(q^{1-n}/bc;q)_kq^k}{(q,q^{1-n}/b,q^{1-n}/c,q^{-n}x,q/x;q)_k}.
\end{align*}
It is easy to verify that
\begin{align*}
&\hskip -3mm G_{n,k}(x,b,c,q)-G_{n,k}(x,b,cq,q)-
\beta_n G_{n-2,k-1}(xq^{-1},bq,cq,q) =0.
\end{align*}
On the other hand, for $c=q^{1-n}/b$, the identity \eqref{eq:32to54}
reduces to the terminating $q$-binomial theorem (see, for example, \cite[(3.3.6)]{Andrews}).

%%%%%%%%%%%%%%%%%%%%%%%%%%%%%%%%%%%%%%%%%%%%%%%%%%%%%%%%%%%%%%%%%%%%%%%%%%%%%%%%%%

\section{Sears' ${}_4\phi_3$ transformations}
Sears \cite{Sears2} derived the following two transformations of ${}_4\phi_3$ series.
The first one is a $q$-analogue of Whipple's formula \cite{Whipple}. It is also one
of the fundamental formulas in the theory of the basic hypergeometric series.
There are several different proofs in the literature. See Andrews and Bowman \cite{AB},
Ismail \cite{Ismail}, Liu \cite{Liu} and Fang \cite{Fang}. Here we shall give
another simple proof of Sears' ${}_4\phi_3$ transformations.

\begin{thm}[Sears' ${}_4\phi_3$ transformations]If $n\geq0 $, there holds
\begin{align}
&\hskip -2mm  _4\phi_3 \left[\begin{array}{c}
q^{-n},\, a,\, b,\, c \\
d,\, e,\, f\end{array}; q,\,q
\right]   \nonumber\\[5pt]
& =\frac{(e/a,f/a;q)_{n}}{(e,f;q)_{n}} a^{n} {}_4\phi_3
\left[\begin{array}{c}
q^{-n}, \, a,\, d/b,\, d/c \\
d,\,aq^{1-n}/e,\, aq^{1-n}/f,\end{array};q,\,q\right]\label{eq:4phi3=01}   \\[5pt]
& =\frac{(a,ef/ab,ef/ac;q)_{n}}{(e,f,ef/abc;q)_{n}}{}_4\phi_3
\left[\begin{array}{c}
q^{-n}, \, e/a,\, f/a,\, ef/abc \\
ef/ab,\,ef/ac,\, q^{1-n}/a,\end{array};q,\,q\right]\label{eq:4phi3=02},
\end{align}
where  $def=abcq^{1-n}$.
\end{thm}

\pf
Let
$$
F_{n,k}(a,b,d,e,f,q)=\frac{(q^{-n},a,b,
def/abq^{1-n};q)_k}{(q,d,e,f;q)_k}q^k
$$
be the $k$-th term in the left-hand side of \eqref{eq:4phi3=01}.
One may check that
\begin{align}
F_{n,k}(aq,b,d,e,f,q)-F_{n,k}(a,b,d,e,f,q)
=\alpha_n
F_{n-1,k-1}(aq,bq,dq,eq,fq,q), \label{eq:fnk4phi3-1}
\end{align}
where
\begin{align*}
\alpha_n &=\frac{(1-q^{-n})(1-b)(a-defq^{n-2}/ab)q}
{(1-d)(1-e)(1-f)}.
\end{align*}

Let
\begin{align*}
G_{n,k}(a,b,d,e,f,q)
=\frac{(e/a,f/a;q)_n a^{n}}{(e,f;q)_n}
\frac{(q^{-n},a,d/b, abq^{1-n}/ef;q)_k}{(q,d,aq^{1-n}/e,aq^{1-n}/f;q)_k}q^k.
\end{align*}
It is easy to verify that
\begin{align}
&\hskip -3mm
G_{n,k}(aq,b,d,e,f,q)-G_{n,k}(a,b,d,e,f,q)-\alpha_n G_{n-1,k-1}(aq,bq,dq,eq,fq,q)  \nonumber\\[5pt]
&=H_{n,k}-H_{n,k-1},  \label{eq:fnk4phi3-2}
\end{align}
where
\begin{align*}
H_{n,k}
&=\frac{a^{n}(a^{2}q^{k-n+2}/ef-1)(e/a,f/a;q)_{n}(q^{-n};q)_{k+1}(aq,d/b,abq^{2-n}/ef;q)_{k}}
{(e,f;q)_{n}(q,d;q)_{k}(aq^{1-n}/e,aq^{1-n}/f;q)_{k+1}}.
\end{align*}
It follows from \eqref{eq:fnk4phi3-1} and \eqref{eq:fnk4phi3-2} that \eqref{eq:4phi3=01}
is true for all $a=q^{-M}$ ($M\geq 0$).

To prove \eqref{eq:4phi3=02}, let
\begin{align*}
P_{n,k}(a,b,d,e,f,q)
=\frac{(a,ef/ab,bq^{1-n}/d;q)_n}{(e,f,q^{1-n}/d;q)_n}
\frac{(q^{-n}, e/a, f/a, q^{1-n}/d;q)_k}{(q, ef/ab, bq^{1-n}/d, q^{1-n}/a;q)_k}
\end{align*}
Then we have
\begin{align}
&\hskip -3mm
P_{n,k}(aq,b,d,e,f,q)-P_{n,k}(a,b,d,e,f,q)-\alpha_n P_{n-1,k-1}(aq,bq,dq,eq,fq,q)  \nonumber\\[5pt]
&=Q_{n,k}-Q_{n,k-1},  \label{eq:fnk4phi3-3}
\end{align}
where
\begin{align*}
Q_{n,k} =
\frac{q^n(efq^{k-1}-a^{2}b)(aq,ef/ab;q)_{n-1}(bq^{1-n}/d;q)_{n}(q^{-n};q)_{k+1}(e/a, f/a, q^{1-n}/d;q)_{k}}
{ab(e,f,q^{1-n}/d;q)_{n}(q, ef/ab, bq^{1-n}/d, q^{1-n}/a;q)_{k}}.
\end{align*}
It follows from \eqref{eq:fnk4phi3-1} and \eqref{eq:fnk4phi3-3} that \eqref{eq:4phi3=02}
is true for all $a=q^{-M}$ ($M\geq 0$).
\qed

%%%%%%%%%%%%%%%%%%%%%%%%%%%%%%%%%%%%%%%%%%%%%%%%%%%%%%%%%%%%%%%%%%%%%%%%%%%%%%%%%%%%%%%%%%%%%%%%%%%%%%%%%%%%%%%%
\section{Chen's bibasic summation formula}
Chu \cite{Chu} established the following bibasic summation formula:
\begin{align}
&\hskip -2mm
\sum_{k=-m}^n\frac{(1-\alpha a_k b_k)(b_k-{a_k}/{\alpha d})}
{(1-\alpha a_0 b_0)(b_0-{a_0}/{\alpha d})}
\frac{\prod_{j=0}^{k-1}(1-a_j)(1-a_j/d)(1-cb_j)(1-\alpha^2db_j/c)}
{\prod_{j=1}^{k}(1-\alpha b_j)(1-\alpha a_j/c)(1-\alpha db_j)(1-ca_j/\alpha d)} \nonumber  \\[5pt]
&=\frac{(1-a_0)(1-a_0/d)(1-cb_0)(1-\alpha^2db_0/c)}
{\alpha(1-\alpha a_0b_0)(1-c/\alpha)(b_0-a_0/\alpha d)(1-\alpha d/c)}  \nonumber \\
&\quad{}\cdot\left(\prod_{j=1}^n
\frac{(1-a_j)(1-a_j/d)(1-cb_j)(1-\alpha^2db_j/c)}
{(1-\alpha b_j)(1-\alpha a_j/c)(1-\alpha db_j)(1-ca_j/\alpha d)}\right.  \nonumber\\
&\qquad{}-\left.
\prod_{j=-m}^{0}
\frac{(1-\alpha b_j)(1-\alpha a_j/c)(1-\alpha db_j)(1-ca_j/\alpha d)}
{(1-a_j)(1-a_j/d)(1-cb_j)(1-\alpha^2 db_j/c)}\right), \label{eq:chu-bi}
\end{align}
which is a generalization of Gasper and Rahman's bibasic summation formula \cite[(3.6.13)]{GR}.

Setting $c=\alpha/B_n$ and $b_k=B_k/\alpha$ for $k=0,1,\ldots,$
Chen \cite{Chen} noticed that \eqref{eq:chu-bi} can be simplified as
\begin{align}
F(n)=\sum_{k=0}^n\frac{(1-a_kB_k)(dB_k-a_k)}{(1-a_0B_n)(dB_n-a_0)}
\frac{\sum_{j=0}^{k-1}(B_n-B_j)(1-dB_n B_j)}{\prod_{j=1}^{k}(1-B_n a_j)(dB_n-a_j)}G(k), \label{eq:chen0}
\end{align}
where
$$
F(n)=\frac{(1-B_0)(1-dB_0)}{(1-B_n)(1-dB_n)},\qquad
G(n)=\frac{\prod_{j=0}^{n-1}(1-a_j)(d-a_j)}{\prod_{j=1}^n(1-B_j)(1-dB_j)}.
$$
Then applying Krattenthaler's matrix inverse \cite{Krattenthaler} to \eqref{eq:chen0},
Chen \cite{Chen} obtained
\begin{align}
G(n)=\sum_{k=0}^n\frac{\prod_{j=0}^{n-1}(1-B_k a_j)(dB_k-a_j)}
{\prod_{j=0,\,j\neq k}(1-dB_k B_j)(B_k-B_j)}F(k). \label{eq:chen2}
\end{align}
In particular, when $a_k=ap^k$ and $B_k=bq^k$, the identity \eqref{eq:chen2}
reduces to the following bibasic summation formula.
\begin{thm}[Chen's bibasic summation formula] If $n\geq0 $, there holds
\begin{align}
\sum_{k=0}^{n}\frac{(1-b)(1-db)(abq^{k},a/dbq^{k};p)_{n}(-1)^{k}q^{{k+1 \choose 2}}}
{(1-bq^{k})(1-dbq^{k})(q,db^{2}q^{k};q)_{k}(q,db^{2}q^{2k+1};q)_{n-k}}
=\frac{(a,a/d;p)_{n}}{(bq,dbq;q)_{n}}.  \label{eq:chenpq}
\end{align}
\end{thm}

\medskip
\pf  Let
$$
F_{n,k}(a,b,d,q)
=\frac{(1-b)(1-db)(abq^{k},a/dbq^{k};p)_{n}(-1)^{k}q^{{k+1 \choose 2}}}
{(1-bq^{k})(1-dbq^{k})(q,db^{2}q^{k};q)_{k}(q,db^{2}q^{2k+1};q)_{n-k}}.
$$
Similarly to the contiguous relation \eqref{eq:contiguous1}, we have
\begin{align}
F_{n,k}(a,b,d,q)-\alpha_n F_{n-1,k}(a,b,d,q)=\beta_n
F_{n-1,k-1}(a,bq,d,q), \label{eq:fnkpq}
\end{align}
where
\begin{align*}
\alpha_n
&=\frac{(1-abp^{n-1})(1-ap^{n-1}/bd)}{(1-q^n)(1-db^2q^n)}, \\
\beta_n
&=\frac{(abq^np^{n-1}-1)(bdq^n-ap^{n-1})(1-bd)(1-b)}{db(1-bdq)(1-bq)(1-db^2q^n)(1-q^n)}.
\end{align*}
The proof then follows from summing \eqref{eq:fnkpq} over $k$ from $0$ to $n$ and proceeding by induction on $n$.
\qed

Note that, the $p=q$ case of \eqref{eq:chenpq} may be rewritten as
\begin{align*}
_8\phi_7 \left[\begin{array}{c}
db^2,\, q\sqrt{db^2},\, -q\sqrt{db^2},\, b ,\,db,\,dbq/a,\,abq^n,\,q^{-n}\\
\sqrt{db^2},\,-\sqrt{db^2},\,
dbq,\,bq,\,ab,\,dbq^{1-n}/a,\,db^2q^{1+n}\end{array}; q,\,q \right]
=\frac{(q,a,a/d,db^2q;q)_{n}}{(ab,bq,dbq,a/db;q)_{n}},
\end{align*}
which is a special case of Jackson's $_{8}\phi_7$ summation formula \cite{Jackson21} (see \cite[(II.22)]{GR}).

%%%%%%%%%%%%%%%%%%%%%%%%%%%%%%%%%%%%%%%%%%%%%%%%%%%%%%%%%%%%%%%%%%%%%%%%%%%%%%%%%%%%
\section{Gasper's split poised $_{10}\phi_9$ transformation}\label{sec:Bailey}
In this section, we show that Gasper's split poised $_{10}\phi_9$
transformation formula (see \cite{Gasper}) can also be proved by the same method.
\begin{thm}[Gasper's split poised $_{10}\phi_9$ transformation] \label{thm:10phi9}
For $n\geq 0$, there holds
\begin{align}
&\hskip -2mm  {}_{10}\phi_9 \left[\begin{array}{c}
a,\,qa^{\frac12},\, -qa^{\frac12},\, b,\, c,\, a/bc,\, C/Aq^n,\, 1/BCq^n,\,B/Aq^n,\,q^{-n} \\
a^{\frac12},\,-a^{\frac12},\,aq/b,\,aq/c,\,bcq,\,1/Cq^n,\,BC/Aq^n,\,1/Bq^n,\,1/Aq^n
\end{array};q,\,q\right]   \nonumber\\[5pt]
&
=\frac{(aq,bq,cq,aq/bc,Aq/B,Aq/C,BCq;q)_n}{(Aq,Bq,Cq,Aq/BC,aq/b,aq/c,bcq;q)_n}
\nonumber\\[5pt]
&\quad{}\times {}_{10}\phi_9 \left[\begin{array}{c}
A,\,qA^{\frac12},\, -qA^{\frac12},\, B,\, C,\, A/BC,\, c/aq^n,\, 1/bcq^n,\,b/aq^n,\,q^{-n} \\
A^{\frac12},\,-A^{\frac12},\,Aq/B,\,Aq/C,\,BCq,\,1/cq^n,\,bc/aq^n,\,1/bq^n,\,1/aq^n
\end{array};q,\,q\right],
\label{eq:10phi9}
\end{align}
where $\lambda=a^2q/bcd$.
\end{thm}

\pf  Let
$$
F_{n,k}(a,b,c,A,B,C,q) =\frac{(a,qa^{\frac12}, -qa^{\frac12}, b, c,
a/bc, C/Aq^n, 1/BCq^n, B/Aq^n,q^{-n};q)_k }
{(a^{\frac12},-a^{\frac12},aq/b,aq/c,bcq,1/Cq^n,BC/Aq^n,1/Bq^n,1/Aq^n;q)_k}q^k.
$$
Similarly to the contiguous relation \eqref{eq:contiguous1} with
$a_r=c$ and $a_{r+1}=a/bcq$, we have
\begin{align}
F_{n,k}(a,b,cq,A,B,C,q)-F_{n,k}(a,b,c,A,B,C,q)=\alpha_n
F_{n-1,k-1}(aq^2,bq,cq,A,B,C,q), \label{eq:10phi9-1}
\end{align}
where
\begin{align*}
\alpha_n &=\frac{(c-a/bcq)(1-bq)(1-aq)(1-aq^2)(1-b)}
{(1-a/c)(1-bcq)(1-aq/b)(1-aq/c)(1-bcq^2)}\\
&\quad\times \frac{(1-C/Aq^n)(1-1/BCq^n)(1-B/Aq^n)(1-q^{-n})q}
{(1-1/Cq^n)(1-BC/Aq^n)(1-1/Bq^n)(1-1/Aq^n)}.
\end{align*}

Let
\begin{align*}
G_{n,k}(a,b,c,A,B,C,q)
&=\frac{(aq,bq,cq,aq/bc,Aq/B,Aq/C,BCq;q)_n}{(Aq,Bq,Cq,Aq/BC,aq/b,aq/c,bcq;q)_n}\\
&\quad{}\times F_{n,k}(A,B,C,a,b,c,q).
\end{align*}
Then we may verify that
\begin{align}
&\hskip -3mm
G_{n,k}(a,b,cq,A,B,C,q)-G_{n,k}(a,b,c,A,B,C,q)-\alpha_n G_{n-1,k-1}(aq^2,bq,cq,A,B,C,q)  \nonumber\\[5pt]
&=H_{n,k}-H_{n,k-1},  \label{eq:hnk00}
\end{align}
where
\begin{align*}
H_{n,k} &=\frac{a(1-b)(1-q^{n-k})(1-aq^{n-k+1})(1-bc^2q/a)}
{b(a-c)(1-aq^{n-k}/bc)(1-bcq^{n-k+1})(1-cq)} G_{n,k}(a,b,c,A,B,C,q).
\end{align*}
Summing \eqref{eq:10phi9-1} and \eqref{eq:hnk00} over $k$ from $0$ to
$n$ respectively, one sees that both sides of \eqref{eq:10phi9}
satisfy the same recurrence relation
\begin{align}
S_{n}(a,b,cq,A,B,C,q)-S_{n}(a,b,c,A,B,C,q)=\alpha_n
S_{n-1}(aq^2,bq,cq,A,B,C,q). \label{eq:rec-cq}
\end{align}
On the other hand, for the case $c=1$, Eq.~\eqref{eq:10phi9} reduces
to
\begin{align*}
1=\frac{(q,Aq/B,Aq/C,BCq;q)_n}{(Aq,Bq,Cq,Aq/BC;q)_n} \sum_{k=0}^n
\frac{(1-Aq^{2k})(A, B, C,
A/BC;q)_k}{(1-A)(q,Aq/B,Aq/C,BCq;q)_k}q^k,
\end{align*}
which is almost trivial since
\begin{align*}
&\hskip -3mm
\frac{(1-Aq^{2k})(A,B, C, A/BC;q)_k}{(1-A)(q,Aq/B,Aq/C,BCq;q)_k}q^k \\
&=\frac{(Aq,Bq,Cq,Aq/BC;q)_k}{(q,Aq/B,Aq/C,BCq;q)_k}
-\frac{(Aq,Bq,Cq,Aq/BC;q)_{k-1}}{(q,Aq/B,Aq/C,BCq;q)_{k-1}}.
\end{align*}
Thus, by \eqref{eq:rec-cq}, Eq.~\eqref{eq:10phi9} holds for all
$c=q^{-M}$ ($M\geq 0$). This completes the proof. \qed

Chu \cite[Corollary 24]{Chu} has generalized Theorem
\ref{thm:10phi9} to a bibasic symmetric transformation as follows.
\begin{thm}[Chu's bibasic symmetric transformation]\label{thm:10phi9chu}
For $n\geq 0$, there holds
\begin{align}
&\hskip -3mm \sum_{k=0}^n\frac{(1-ap^{2k})(a, b, c,
a/bc;p)_k(C/Aq^n, 1/BCq^n, B/Aq^n,q^{-n};q)_k}
{(1-a)(p,ap/b,ap/c,bcp;p)_k (1/Cq^n,BC/Aq^n,1/Bq^n,1/Aq^n;q)_k}p^k  \nonumber\\
&=\frac{(ap,bp,cp,ap/bc;p)_n(q,Aq/B,Aq/C,BCq;q)_n}{(p,ap/b,ap/c,bcp;p)_n(Aq,Bq,Cq,Aq/BC;q)_n} \nonumber\\
&\quad\times \sum_{k=0}^n\frac{(1-Aq^{2k})(A, B, C,
A/BC;q)_k(c/ap^n, 1/bcp^n, b/ap^n,p^{-n};p)_k}
{(1-A)(q,Aq/B,Aq/C,BCq;q)_k (1/cp^n,bc/ap^n,1/bp^n,1/ap^n;p)_k}q^k.
\label{eq:chu-10phiq}
\end{align}
\end{thm}
We point out that \eqref{eq:chu-10phiq} can be proved in the same
way as \eqref{eq:10phi9}. Let $F_{n,k}(a,b,c,A,B,C;p,q)$ and
$G_{n,k}(a,b,c,A,B,C;p,q)$ be the $k$-th summands in the left-hand
side and right-hand side of \eqref{eq:chu-10phiq} respectively. Then
there exist relations for $F_{n,k}$ and $G_{n,k}$ exactly similar to
\eqref{eq:10phi9-1} and \eqref{eq:hnk00}. For instance, we have
\begin{align*}
F_{n,k}(a,b,cp,A,B,C;p,q)-F_{n,k}(a,b,c,A,B,C;p,q)=\alpha_n
F_{n-1,k-1}(ap^2,bp,cp,A,B,C;p,q),
\end{align*}
where
\begin{align*}
\alpha_n &=\frac{(c-a/bcp)(1-bp)(1-ap)(1-ap^2)(1-b)}
{(1-a/c)(1-bcp)(1-ap/b)(1-ap/c)(1-bcp^2)}\\
&\quad\times \frac{(1-C/Aq^n)(1-1/BCq^n)(1-B/Aq^n)(1-q^{-n})p}
{(1-1/Cq^n)(1-BC/Aq^n)(1-1/Bq^n)(1-1/Aq^n)}.
\end{align*}

\section{A finite form of Sylvester's identity and Jacobi's triple product identity}
From now on, we assume that $|q|<1$. In the same vein of \eqref{eq:fnk-ind}, we will prove a finite form of Sylvester's identity
and Jacobi's triple product identity.
Recall that the {\it $q$-binomial coefficients} are defined by
\begin{align*}
{n\brack k}=\frac{(q;q)_n}{(q;q)_k(q;q)_{n-k}},\quad n,k\in\mathbb{Z}.
\end{align*}

\begin{thm}[A finite form of Sylvester's and Jacobi's identities]
For $n\geq 0$, there holds
\begin{align}
\sum_{k=0}^{n}(-1)^k x^kq^{k(3k+1)/2}{n\brack k}
\frac{(1-x q^{2k+1})}{(xq^{k+1};q)_{n+1}}=1. \label{eq:sylvester-1}
\end{align}
\end{thm}
\pf Let
\begin{align*}
F_{n,k}(x,q)=(-1)^k x^kq^{k(3k+1)/2}{n\brack k}
\frac{(1-x q^{2k+1})}{(xq^{k+1};q)_{n+1}}.
\end{align*}
Noticing the trivial relation
\begin{align}
{n\brack k}(1-xq^{n+1})
={n-1\brack k}(1-xq^{n+k+1})+{n-1\brack k-1}(1-xq^{k+1})q^{n-k},
\label{eq:bino-rec}
\end{align}
we have
\begin{align*}
F_{n,k}(x,q)=\frac{1}{1-xq^{n+1}}F_{n-1,k}(x,q)
-\frac{xq^{n+1}}{1-xq^{n+1}}F_{n-1,k-1}(xq^2,q).
\tag*{\qed}
\end{align*}

Letting $n\rightarrow\infty$ in
\eqref{eq:sylvester-1}, we immediately obtain Sylvester's identity
\cite[(9.2.3)]{Andrews}:
\begin{align}
\sum_{k=0}^{\infty} \frac{(-1)^k x^kq^{k(3k+1)/2}(1-x
q^{2k+1})}{(q;q)_k(xq^{k+1};q)_{\infty}}=1.  \label{eq:sylvester}
\end{align}
On the other hand, if we make the substitutions $n\rightarrow m+n$,
$x\rightarrow -xq^{-2m}$ and $k\rightarrow m+k$ in
\eqref{eq:sylvester-1}, and apply the following relation
\begin{align}
\frac{(-xq^{-2m})^{m+k}
q^{(3(m+k)^2+m+k)/2}}{(xq^{k-m+1};q)_{m+n+1}}
&=\frac{x^{2k}q^{2k^2+k}}{(1/x;q)_{m-k}(xq;q)_{n+k+1}},
\label{eq:relation}
\end{align}
we obtain
\begin{align}
&\hskip -3mm \sum_{k=-m}^{n}{m+n\brack m+k}
\frac{(1-xq^{2k+1})x^{2k}q^{2k^2+k}}{(1/x;q)_{m-k}(xq;q)_{n+k+1}}=1.
\label{eq:syl-2}
\end{align}
Letting $m,n\rightarrow \infty$ in  \eqref{eq:syl-2}, we immediately
get
\begin{align}
\sum_{k=-\infty}^\infty
(1-xq^{2k+1})x^{2k}q^{2k^2+k}=(1/x,xq,q;q)_\infty,
\label{eq:jacobi2k}
\end{align}
which is Jacobi's triple product identity (see \cite[p.~21]{Andrews}).

Note that finite forms of the quintuple product identity (see \cite[p.~147]{GR})
\begin{align*}
\sum_{k=-\infty}^{\infty}(z^2 q^{2k+1}-1)
z^{3k+1}q^{k(3k+1)/2}=(q,z,q/z;q)_\infty (qz^2,q/z^2;q^2)_\infty,
\end{align*}
which is related to Jacobi's triple product identity, were given by \cite{Paule,CCG,GZ07}.

\section{A finite form of Kang's identity}
In this section, we give a finite form of Kang's \cite{Kang} generalization of
Sylvester's identity.  Similarly to the previous section, we only give the corresponding contiguous relation
for the $k$th terms in our proof.

\begin{thm}For $n\geq 0$, there holds
\begin{align}
\sum_{k=0}^{n}(-1)^k x^kq^{k(3k+1)/2}{n\brack k}
\frac{(cq^{-k},dq^{-k};q)_k(1-xq^{2k+1})}{(cx,dx;q)_{k+1}(xq^{k+1};q)_{n+1}}
=\frac{(cdx;q)_n}{(cx,dx;q)_{n+1}}. \label{eq:sylvester-2}
\end{align}
\end{thm}
\pf Let
\begin{align*}
F_{n,k}(c,d,x,q)=(-1)^k x^kq^{k(3k+1)/2}{n\brack k}
\frac{(cq^{-k},dq^{-k};q)_k(1-xq^{2k+1})}{(cx,dx;q)_{k+1}(xq^{k+1};q)_{n+1}}.
\end{align*}
Then
\begin{align*}
F_{n,k}(c,d,x,q) &=\frac{1}{1-xq^{n+1}}F_{n-1,k}(c,d,x,q)  \\
&\qquad{}-\frac{(1-cq^{-1})(1-dq^{-1})xq^{n+1}}{(1-cx)(1-dx)(1-xq^{n+1})}F_{n-1,k-1}(cq^{-1},dq^{-1},xq^2,q).
\end{align*}
\qed

Letting $n\rightarrow\infty$ in \eqref{eq:sylvester-2}, we obtain
the following identity due to Kang \cite[(6.9)]{Kang}: %which is a generalization of Sylvester's identity \eqref{eq:sylvester}:
\begin{align*}
\sum_{k=0}^{\infty}(-1)^k x^k q^{k(3k+1)/2}
\frac{(aq^{-k},bq^{-k},xq;q)_k(1-xq^{2k+1})}{(q;q)_k(ax,bx;q)_{k+1}}
=\frac{(abx,xq;q)_\infty}{(ax,bx;q)_\infty}
\end{align*}
for $|ax|<1$ and $|bx|<1$.

%\vskip 5mm \noindent{\bf Acknowledgement.} We thank the referee for helpful comments on a previous
%version of this paper.

\renewcommand{\baselinestretch}{1}

\end{document}